\documentclass[a4paper]{article}

\usepackage[breaklinks, colorlinks, linkcolor=blue, citecolor=blue, urlcolor=black]{hyperref}

\usepackage{graphicx}
\usepackage{todonotes}
\usepackage{tabularx}
\usepackage{tabularcalc}
\usepackage{booktabs}
\usepackage{multirow}
\usepackage{subfig}
\usepackage{siunitx}
\usepackage{soul}
\usepackage{footmisc}
\usepackage{svg}
\usepackage{listings}
\usepackage{algorithm}
\usepackage{algorithmic}
\usepackage{dsfont}
\usetikzlibrary{calc,positioning}

\usepackage{zibtitlepage}
\ZTPTitle{Combining Precision Boosting with LP Iterative Refinement for Exact Linear Optimization}
\ZTPAuthor{ \ZTPHasOrcid{Leon Eifler}{0000-0003-0245-9344} Jules Nicolas-Thouvenin \ZTPHasOrcid{Ambros Gleixner}{0000-0003-0391-5903}
}
\ZTPNumber{ZIB-Report 23-26}
\ZTPMonth{November}
\ZTPYear{2023}
\usepackage[margin=1.5in]{geometry}

\svgpath{{pictures/}} 
\usepackage{graphicx}
\newcommand{\myorcidlink}[1]{\,\href{https://orcid.org/#1}{\raisebox{-0.45ex}{\includegraphics[width=1.8ex]{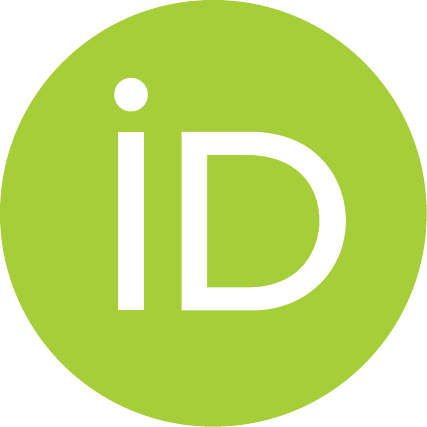}}}}

\usepackage{amsmath}
\usepackage{amssymb}
\usepackage{amsthm}
\usepackage{xspace}
\usepackage{xcolor}
\usepackage{collcell}
\usepackage{accents}

\newcommand{\ie}{i.e.,\xspace}
\newcommand{\eg}{e.g.,\xspace}

\newcommand{\Abs}[1]{\vert #1\vert}
\newcommand{\Norm}[1]{\Vert #1\Vert}

\newcommand{\floor}[1]{\lfloor #1\rfloor}

\newcommand{\onevec}{\mathds{1}}

\newcolumntype{R}{>{\collectcell\ApplyColor}{r}<{\endcollectcell}}
\newcommand{\solver}[1]{\textsc{#1}\xspace}

\newcommand{\scipversion}{8.0}

\newcommand{\scip}{\solver{SCIP}}
\newcommand{\scipv}{\solver{SCIP}~\scipversion\xspace}
\newcommand{\soplex}{\solver{SoPlex}}
\newcommand{\soplexversion}{6.0.3}
\newcommand{\soplexv}{\solver{SoPlex}~\soplexversion\xspace}

\newcommand{\papiloversion}{2.0.1}
\newcommand{\papilov}{\solver{PaPILO}~\papiloversion\xspace}

\newcommand{\qsoptex}{\solver{QSopt\_ex}}

\newcommand{\miplib}{\mbox{MIPLIB}}

\newcommand{\safegmi}{\textsc{safegmi+ir-double}\xspace}
\newcommand{\safegmiboost}{\textsc{safegmi+ir-boosting}\xspace}
\newcommand{\enc}{\textsc{densmall+ir-double}\xspace}
\newcommand{\encboost}{\textsc{densmall+ir-boosting}\xspace}

\newcommand{\combi}{\textsc{ir-boosting}\xspace}
\newcommand{\pureir}{\textsc{ir-double}\xspace}
\newcommand{\purepb}{\textsc{boosting-pure}\xspace}

\newcommand{\irpaper}{\textsc{lplib}\xspace}
\newcommand{\cutlib}{\textsc{cutlib}\xspace}

\newcommand{\Q}{\mathbb{Q}\xspace}

\newcommand{\percentage}[2]{\pgfmathparse{#1/#2 *100}\pgfmathprintnumber[precision=1]{\pgfmathresult}\%}
\newcommand{\reduction}[2]{\pgfmathparse{(#2-#1)/#2 *100}\pgfmathprintnumber[precision=1]{\pgfmathresult}\%}
\newcommand{\increase}[2]{\pgfmathparse{(#1-#2)/#2 *100}\pgfmathprintnumber[precision=1]{\pgfmathresult}\%}
\newcommand{\fraction}[2]{\pgfmathparse{(#1)/#2}\pgfmathprintnumber[precision=1]{\pgfmathresult}}

\newcommand{\colorquot}[2]{
   \pgfmathparse{(#1<0.9*#2)?1:0}\ifdim\pgfmathresult pt>0pt \textcolor{blue}{$\mathbf{\reduction{#1}{#2}}$}\else
      \pgfmathparse{(#1>1.1*#2)?1:0}\ifdim\pgfmathresult pt>0pt \textcolor{red}{$\mathbf{\reduction{#1}{#2}}$} \else
         $\reduction{#1}{#2}$\fi \fi
}

\newlength\myindent
\newcommand\bindent{%
   \begingroup
   \setlength{\itemindent}{\myindent}
   \addtolength{\algorithmicindent}{\myindent}
}
\newcommand\eindent{\endgroup}

\newcommand{\K}{10000}

\newcommand{\fp}{floating-point\xspace}

\newcommand{\pb}{precision boosting\xspace}

\newcommand{\ir}{LP iterative refinement\xspace}

\begin{document}


\title{Combining Precision Boosting with LP Iterative Refinement for Exact Linear Optimization\thanks{The work for this article has been conducted within the Research Campus Modal funded by the German Federal Ministry of Education and Research (BMBF grant numbers 05M14ZAM, 05M20ZBM).}}
\author{Leon Eifler$^1$\myorcidlink{0000-0003-0245-9344} \and Jules Nicolas-Thouvenin \and
   Ambros Gleixner$^{1,2}$\myorcidlink{0000-0003-0391-5903}
}
\date{%
   $^1$Zuse Institute Berlin\\%
   $^2$HTW Berlin\\[2ex]%
   $^3$
   \today
}

\maketitle

\begin{abstract}
   This article studies a combination of the two state-of-the-art algorithms for the exact solution of linear programs (LPs) over the rational numbers, \ie without any roundoff errors or numerical tolerances.
   By integrating the method of \emph{precision boosting} inside an \emph{\ir} loop, the combined algorithm is able to leverage the strengths of both methods: the speed of \ir, in particular in the majority of cases when a double-precision floating-point solver is able to compute approximate solutions with small errors, and the robustness of precision boosting whenever extended levels of precision become necessary.
   We compare the practical performance of the resulting algorithm with both pure methods on a large set of LPs and mixed-integer programs (MIPs). The results show that the combined algorithm solves more instances than a pure \ir approach, while being faster than pure precision boosting. When embedded in an exact branch-and-cut framework for MIPs, the combined algorithm is able to reduce the number of failed calls to the exact LP solver to zero, while maintaining the speed of the pure \ir approach.
\end{abstract}

\section{Introduction}
\label{sec:intro}

Linear programming (LP) is a fundamental optimization technique widely used in various fields, including operations research, engineering, economics, and finance.
In practice, linear programming solvers rely on fast floating-point arithmetic, coupled with the careful use of error tolerances to efficiently compute accurate solutions.
However, the use of floating-point arithmetic can lead to numerical inaccuracies, especially for problems with large coefficient ranges, which in turn can result in inaccurate solutions or incorrect claims of optimality or infeasibility. Exact linear programming algorithms aim to solve LPs exactly over the rational numbers, \ie without any numerical inaccuracies or error tolerances. Such exact solvers are needed as a subroutine for exact mixed integer programming (MIP) \cite{CookKochSteffyetal2013, EiflerGleixner2022}, but can also directly be used to investigate numerically challenging LPs or to establish theoretical results \cite{Hales2017, Lerman2012, Bofi19, Burt12, EiflerGleixnerPulaj2022, KenterEtAl2018, LanciaEtAl2020,Pulaj20}.

The na\"ive approach of performing a simplex method in exact arithmetic was observed to be prohibitively slow in many practical applications by Espinoza \cite{Espinoza2006}.
Also the idea of using limited-precision arithmetic at a fixed, but sufficiently high level in order to obtain theoretical guarantees of convergence to an exact solution is limited in its practical applicability.
This holds even for algorithms with polynomial runtime, such as the algorithm described by Gr\"otschel, Lov\'asz, and Schrijver~\cite{GroetschelLovaszSchrijver1988}, which is itself based on the ellipsoid method of Khachiyan~\cite{KHACHIYAN198053}, see~\cite{GleixnerSteffy2020}.
In this article, we focus purely on simplex-based methods.
Their warm-starting capabilities also align well with one of over main motivations to us exact LP solvers as a subroutine in exact branch-and-bound solvers. 

Among simplex-based methods for solving LPs exactly, the more successful approaches rely on combining floating-point arithmetic and exact arithmetic in some way. An early attempt was presented by Dhiflaoui et al.~\cite{Dhiflaoui2003}. They first solve the LP approximately in double-precision arithmetic and try to prove the optimality of the found solution by symbolically factorizing the returned basis matrix. If this approach fails, they continue with an exact rational simplex, warm-started from the final floating-point basis.
This approach was refined by Applegate et al \cite{Espinoza2006, APPLEGATE2007} in an algorithm called incremental precision boosting, which is implemented in the solver \qsoptex. In each iteration of the algorithm, the floating-point precision is increased until the basis can be proven to be exactly optimal.
Another state-of-the-art algorithm for solving LPs exactly is based on \ir \cite{GleixnerSteffyWolter2016, GleixnerSteffy2020} and is implemented in the LP solver \soplex. It avoids higher-precision LP solves by instead solving a series of error-correcting LPs in double-precision in order to produce a sequence of primal-dual solutions with residual errors converging to zero.

Although \ir was shown to outperform precision boosting in the majority of cases, the precision boosting algorithm is more robust on numerically difficult instances. The \ir procedure has no reliable way of recovering if the floating-point LP subroutine aborts with a failure due to numerical issues. This has been observed in practice both on pure LPs~\cite{GleixnerSteffyWolter2016,GleixnerSteffy2020} and for LP relaxations of exact MIP subproblems during branch-and-cut~\cite{EiflerGleixner2023}.

In this article, we propose a natural combination of these two algorithms that profits from the speed of \ir, but can use precision boosting as a fallback to overcome numerical issues. We show that this combination is more robust and faster than either of the algorithms individually.

The remainder of this paper is structured as follows. In Section~\ref{sec:review} we introduce exact LP solving formally and give a brief review of the two base algorithms. In Section~\ref{sec:contribution}, we present the combined algorithm, giving details on when and how the precision boosting technique is used. In Section~\ref{sec:comp}, we conduct a computational study, evaluating the different algorithms both in the context of pure exact LP, as well as in experiments with an exact MIP framework.
Finally, in Section~\ref{sec:conc}, we conclude our findings and give an outlook on future work.


\section{Existing Methods to Solve LPs Exactly}
\label{sec:review}

We aim to solve a linear program
\begin{align}
   \label{eq:lpprimal}
   \min \{c^Tx | Ax = b, x \geq \ell\}
\end{align}
where $A \in \Q^{m \times n}$ is a rational matrix of full row rank with $m \le n$,
$c \in \Q^n$ is the objective function, and $\ell \in \Q^m$ is the lower bound vector. Note that we choose this formulation to keep the notation simple.  More general formulations are possible and discussed in detail, e.g., in \cite{GleixnerSteffyWolter2016}.

Our goal is to solve this LP exactly over the rational numbers, \ie any feasibility or optimality tolerances that are often used in LP solvers are set to zero. As discussed in Section~\ref{sec:intro}, we focus on algorithms based on the simplex method. Both LP iterative refinement and precision boosting iteratively produce a sequence $(x_n,y_n)$ of approximate solutions that are more and more accurate.
With any such algorithm, we can then use the methods described in \cite{GleixnerSteffy2020} to obtain an exact solution, either by reconstructing it from an approximate solution as described in \cite{GleixnerSteffy2020}, or by solving the linear system defining the current basis exactly.

In the following, we describe the two existing methods that we combine in this paper. We refer to \cite{Espinoza2006} for a more detailed description of precision boosting and to \cite{GleixnerSteffyWolter2016, GleixnerSteffy2020} for a more detailed description of \ir.


\subsection{Incremental precision boosting}
This algorithm first computes an approximate solution using a floating-point simplex implementation in double-precision arithmetic and checks the resulting basis for exact primal and dual feasibility.
In the original algorithm \cite{Espinoza2006} this check is always performed by means of a rational LU factorization of the basis matrix and subsequent triangular solves in rational arithmetic.
If the basis is detected as not optimal, the arithmetic precision is increased and all floating-point tolerances are decreased. After double precision, the first level of extended precision uses 128~bits (quad precision); subsequently, the precision is grown by a factor of $1.5$ in each boosting step.

The tolerances are decreased at the same rate as the precision is increased.
Given a value of $2^a$ for some $a \in \Q$ as the tolerance's default value in double precision and $p$ bits for the mantissa in higher precision, the corresponding tolerance value is set to $2^{a\frac{p}{64}}$. Note that this leaves a "buffer", as the mantissa in double precision only has $53$ bits. As an example, a tolerance value of $10^{-6}$ in double precision would be scaled to approximately $3\cdot 10^{-11}$ in quad precision.

If an iteration returns that the LP is infeasible, the algorithm attempts to turn the approximate Farkas proof into an exact proof of infeasibility.
A short algorithmic description for a feasible LP can be found in Algorithm \ref{alg:boosting}.

\begin{algorithm}
   \caption{Incremental precision boosting}
   \label{alg:boosting}
   \begin{algorithmic}
      \STATE{\textbf{Input}: $c,\ell \in \Q^n,b \in \Q^m, A\in \Q^{n \times m}$ }
      \STATE{\textbf{Output}: primal-dual solution $(x^*,y^*) \in \Q^{n+m}$ of \eqref{eq:lpprimal}, basis
         $\mathcal{B}$}
      \bindent
      \FOR{$p \gets 64,128,192,288,\ldots$}
      \STATE{load $\bar A,\bar b,\bar c,\bar\ell$ with precision $p$}
      \STATE{decrease floating-point tolerances}
      \STATE{solve $\min \{\bar c^Tx | \bar Ax = \bar b, x \geq \bar\ell\}$ in precision $p$}
      \STATE{$\mathcal{B} \gets$ returned basis}
      \STATE{symbolically compute solution $(x^*,y^*)$  corresponding to $\mathcal{B}$}
      \IF{$(x^*,y^*)$ is primal and dual feasible}
      \RETURN{$(x^*,y^*),\mathcal{B}$}
      \ENDIF
      \ENDFOR
      \eindent
   \end{algorithmic}
\end{algorithm}

The \pb algorithm is often very efficient, as it has been observed that the basis returned by the first floating-point LP solve is often already exactly optimal \cite{Dhiflaoui2003}. Furthermore, the higher numeric precision on successive iterates makes this algorithm very robust on problems that are numerically difficult, even in the presence of large coefficient ranges or ill-conditioned basis matrices.
Its downside is that the computations in higher precision can be time-consuming, and the exact factorization of the basis matrix after every iteration can pose an additional bottleneck in some problem instances.

\subsection{\ir}
In an attempt to overcome these issues, the \emph{\ir} algorithm was introduced \cite{GleixnerSteffyWolter2016}, which is based on iterative refinement for linear systems \cite{Wilkinson1963}.
Instead of increasing the numerical precision of the floating-point computations, \ir computes the violations in the reduced costs, right-hand sides, and the lower bounds. Then, those residuals are scaled and inserted in place of the original objective, right-hand sides, and lower bounds, respectively.

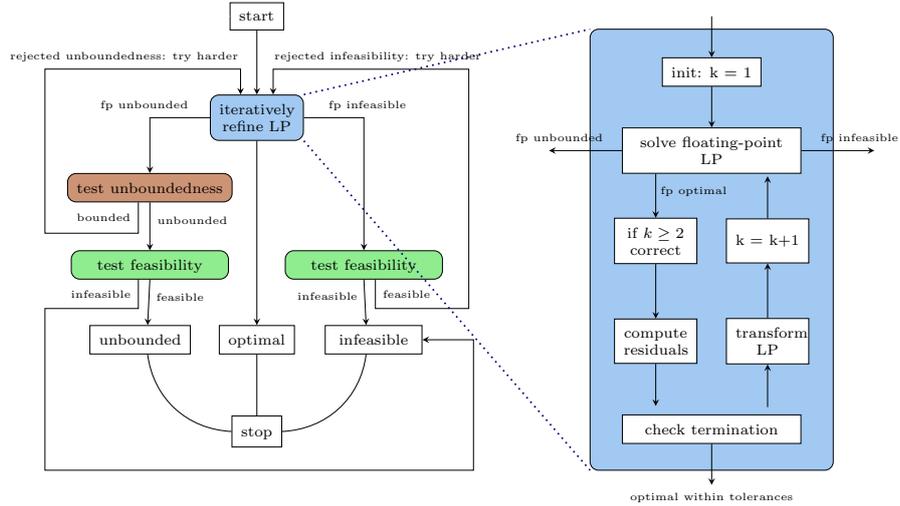
\begin{figure}[h]
   \caption{Iterative refinement for exact linear programming, reprinted from \cite{GleixnerSteffyWolter2016}}
   \label{fig:ir}
   \centering
   \scalebox{0.8}{
\definecolor{DaLightBlue}{rgb}{0.63,0.79,0.95}
\definecolor{DaLightBrown}{rgb}{0.8,0.58,0.46}
\definecolor{DaLightGreen}{rgb}{0.56,0.93,0.56}
\definecolor{DaDarkBlue}{rgb}{0.0,0.0,0.46}
\begin{tikzpicture}[auto, scale=.8, >=stealth]
  \tikzstyle{state} += [draw=black,rectangle, inner sep=1.5mm];
  \tikzstyle{link} += [->];
  \node (start) at (5.5,6.3) [state] {\scriptsize start};
  \node (refine) at (5.5,4.2) [state,fill=DaLightBlue,rounded corners] {\parbox{3.5em}{\scriptsize\centering iteratively\\refine LP}};
  \node (unb) at (3.1,-0.4) [state] {\scriptsize unbounded};
  \node (opt) at (5.5,-0.4) [state] {\scriptsize \smash[b]{optimal}};
  \node (inf) at (7.9,-0.4) [state, inner xsep=1.5ex] {\scriptsize infeasible};
  \node (stop) at (5.5,-2.3) [state, inner xsep=1ex, inner ysep=1.2ex] {\scriptsize \smash[b]{stop}};
  \node (testunbd) at (3.3,2.75) [state, rounded corners,fill=DaLightBrown] {\scriptsize test unboundedness};
  \node (testfeasunbd) at (3.3,1.15) [state, inner xsep=2.8ex,rounded corners,fill=DaLightGreen] {\scriptsize \smash[b]{test feasibility}};
  \node (testfeas) at (7.7,1.15) [fill=white, state, inner xsep=2.8ex, rounded corners,fill=DaLightGreen] {\scriptsize \smash[b]{test feasibility}};

  \draw [link] (start) -- (refine);
  \draw [link] (refine) -- (opt); 
  \draw (unb.-63) edge [bend right=40] (stop.west);
  \draw (opt) -- (stop);
  \draw (inf.-117) edge [bend left=40] (stop.east);

  \draw[white] (start) -- ++(-5.5cm,0); 

  \draw [link] (refine.west) -| node[anchor=south] {\tiny\hspace*{-3.5ex} fp unbounded} (testunbd.north);
  \draw [link] (testunbd.south) -- node[anchor=188] {\tiny unbounded} (testfeasunbd.north);
  \draw [link] (testunbd.232) -- node[anchor=east] {\tiny bounded} ++(0,-0.64) -| ++(-1.92,3.4)
    -| node[anchor=south, near start] {\tiny\hspace*{-2.9em} rejected unboundedness: try harder} (refine.125);
    \draw [link] (testfeasunbd.south) -- node[anchor=190.7] {\tiny feasible} (unb.63);
  \draw [link] (testfeasunbd.232) -- node[anchor=east] {\tiny infeasible\hspace*{0ex}} ++(0,-0.6) -| ++(-1.92,-3.35) -- ++(8.81,0) |- (inf.east);

  \draw [link] (refine.east) -| node[anchor=south] {\tiny fp infeasible\hspace*{-1.5ex}} (testfeas.north);
  \draw [link] (testfeas.308) -- node[anchor=west] {\tiny feasible\hspace*{0ex}} ++(0,-0.6) -| ++(1.92,4.958)
    -| node[anchor=south, near start] {\tiny rejected infeasibility: try harder\hspace*{-0.95em}} (refine.55);
  \draw [link] (testfeas.south) -- node[anchor=-9] {\tiny infeasible} (inf.117);

  \begin{scope}[xshift=3.35cm, yshift=-2.86cm, yscale=1.215]
  \node (b1) at (9,7.32) []{};
  \node (b2) at (9,-0.2) []{};
  \draw [black,fill=DaLightBlue, rounded corners] (b1) rectangle (14,-0.2);
  \node (firstk) at (11.5,6.59) [state,fill=white] {\scriptsize init: k = 1};
  \node (fpsolve) at (11.5,5.25) [state,fill=white] {\parbox{7.5em}{\scriptsize\centering solve floating-point LP}};
  \node (correct) at (10.35,3.71) [state,fill=white] {\parbox{3em}{\scriptsize\centering if $k \geq 2$\\correct}};
  \node (residuals) at (10.35,2) [state,fill=white, inner ysep=1ex] {\parbox{3em}{\scriptsize\centering compute\\ residuals}};
  \node (terminate) at (11.5,0.5) [state,fill=white] {\parbox{7.5em}{\scriptsize\centering check termination}};
  \node (trans) at (12.65,2) [state,fill=white] {\parbox{3em}{\scriptsize\centering transform\\LP}};
  \node (nextk) at (12.65,3.71) [state,fill=white, inner ysep=1.75ex] {\parbox{3em}{\scriptsize\centering k = k+1}};

  \draw [link] (firstk) -- (fpsolve);
  \draw [link] (correct.90) ++(0,0.75) -- node[anchor=187] {\tiny \hspace*{-0.6ex}fp optimal} (correct.90);
  \draw [link] (correct) -- (residuals);
  \draw [link] (residuals.270) -- ++(0,-0.72);
  \draw [link] (trans.270) ++(0,-0.73) -- (trans.270);
  \draw [link] (trans) -- (nextk);
  \draw [link] (nextk.90) -- ++(0,0.72);
  \draw [link] (fpsolve.west) -- node {} +(-1.5,0) node[left, anchor=-130] {\tiny fp unbounded};
  \draw [link] (fpsolve.east) -- node[swap] {} +(1.5,0) node[right, anchor=-40] {\tiny fp infeasible};
  \draw [link] (terminate.south) -- node {} +(0,-0.7) node[below, anchor=north] {\tiny optimal within tolerances};
  \draw [<-] (firstk.north) -- +(0,0.7) node[above] {};
\end{scope}

  \draw [dotted,line width=0.8pt,DaDarkBlue] (refine.north east) -- (b1.center);
  \draw [dotted,line width=0.8pt,DaDarkBlue] (refine.south east) -- (b2.center);
\end{tikzpicture}
}
\end{figure}

At each iteration, this transformed problem is solved in double-precision floating-point arithmetic.
Afterwards, the original solution is updated by adding an unscaled version of the transformed solution. This correction is performed in rational arithmetic, and the resulting solution is then used as a starting point for the next iteration.
The algorithm for a feasible LP is provided in Algorithm \ref{fig:ir}.

\begin{algorithm}[h]
   \caption{Iterative refinement for linear programming}
   \label{alg:ir}
   \begin{algorithmic}
      \STATE{\textbf{Input}: $c,\ell \in \Q^n,b \in \Q^m, A\in \Q^{n \times m}$, scaling limit $\alpha$ }
      \STATE{\textbf{Output}: primal-dual solution $(x^*,y^*) \in \Q^{n+m}$ of \eqref{eq:lpprimal}, basis
         $\mathcal{B}$}
      \bindent
      \STATE{$\delta^1_P \gets 1, \delta^1_D \gets 1$}
      \STATE{load $\bar A,\bar b,\bar c,\bar\ell$ in double precision \hfill \COMMENT{initial solve}}
      \STATE{solve $\min \{\bar c^Tx | \bar Ax = \bar b, x \geq \bar\ell\}$ approximately}
      \STATE{$(x^1,y^1) \gets$ returned approximate solution, $\mathcal{B} \gets$ returned Basis}
      \FOR{k $\gets 1,2,\ldots$}
      \STATE{symbolically compute solution $(x^*,y^*)$  corresponding to $\mathcal{B}$ \hfill \COMMENT{check termination}}
      \STATE{or reconstruct $(x^*,y^*)$ from $(x^k,y^k)$}
      \IF{$(x^*,y^*)$ is primal and dual feasible}
      \RETURN{$(x^*,y^*),\mathcal{B}$}
      \ELSE
      \STATE{$\hat b \gets b - Ax_k, \hat \ell \gets \ell - x_k$ \hfill \COMMENT{compute error}}
      \STATE{$\hat c \gets c - y_k^T A$}
      \STATE{$\delta^k_P \gets \max\{\Norm{\hat b}_{\infty},\Norm{\hat \ell}_{\infty}\}$}
      \STATE{$\delta^k_D \gets \max\{0, \max\{-\hat c_i | i = 1,\ldots,n\}\}$}
      \STATE{$\Delta^{k+1}_P \gets 1 / \max\{\delta^k_P, (\alpha \Delta^k_P)^{-1}\}$ \hfill \COMMENT{compute scaling}}
      \STATE{$\Delta^{k+1}_D \gets 1 / \max\{\delta^k_D, (\alpha \Delta^k_D)^{-1}\}$}
      \STATE{$\bar b = \Delta^{k+1}_P \hat b, \bar \ell = \Delta^{k+1}_P \hat \ell, \bar c = \Delta^{k+1}_D \hat c$ \hfill \COMMENT{solve transformed}}
      \STATE{solve $\min \{\bar c^Tx | \bar Ax = \bar b, x \geq \bar\ell\}$ approximately }
      \STATE{$(\hat x ,\hat y) \gets $ returned approximate solution}
      \STATE{$x_{k+1} \gets x_k + 1/\Delta^{k+1}_P \hat x$ \hfill \COMMENT{update solution}}
      \STATE{$y_{k+1} \gets y_k + 1/\Delta^{k+1}_D \hat y$}
      \ENDIF
      \ENDFOR
      \eindent
   \end{algorithmic}
\end{algorithm}

In the case that an iteration detects floating-point infeasibility, the \ir algorithm solves the auxiliary \emph{feasibility problem} \cite{GleixnerSteffyWolter2016}
\begin{align}
   \label{eq:feasprob}
   \max\{ \tau\, |\, A\xi - (b-A\ell) \tau = 0,\, \xi \ge 0,\, \tau \le 1\}.
\end{align}
exactly with the previosuly described \ir algorithm. If unboundedness is detected during an iteration, both a primal feasible solution, as well as an unbounded direction of improvement need to be computed. A primal feasible solution can be computed by solving the feasibility problem \eqref{eq:feasprob} as described above. To compute an unbounded ray, the \emph{unboundedness problem}
\begin{align}
   \label{eq:unbdprob}
   Av = 0,\, c^Tv = -1,\, v \ge 0.
   \end{align}
is solved.
A flowchart describing the full algorithm can be seen in Figure \ref{fig:ir}.

The strength of \ir is that it can perform all simplex solves in fast double-precision arithmetic and requires symbolic computations only to recover exact solutions and to create the transformed problems. This results in faster running times on instances that can be solved by \ir\cite{GleixnerSteffyWolter2016}.
The downside is that it cannot reliably recover when the floating-point solver fails due to numerical difficulties. Furthermore, if the floating-point solver reports infeasibility but the feasibility problem \eqref{eq:feasprob} is feasible, there is also no reliable recovery mechanism.
In those cases, the pure \ir algorithm attempts to overcome the difficulties by trying different setting combinations of presolving, scaling, ratio testing, pricing, increasing the Markowitz threshold, as well as changing tolerances.

In order to harness the individual strengths of both algorithms, we propose a combination, described in the following section.

\section{Combining Precision Boosting and \ir}
\label{sec:contribution}

As the \ir algorithm tends to be faster whenever it succeeds, we use it as the main algorithm in our approach, and use precision boosting only whenever \ir fails.
Concretely, this can happen in three cases.

\subsection{Possible failures of pure \ir}

\paragraph{Failure due to numerical troubles.} The first case occurs when the \fp solver, executed inside \ir, fails to terminate with an approximately optimal solution.
Reasons for this can be that a linear system $Bx = b$ needs to be solved where $B$ is numerically singular, or that cycling occurs, \ie the same few variables enter/leave the basis repeatedly. In these cases, the floating-point solver will return correctly with a non-optimal solution status. A different case of numerical difficulties arises when the \fp solver claims to have reached an optimal basis, but in fact returns a solution with large residual errors, \ie close to or even exceeding one. Then, the \ir algorithm can not be guaranteed to converge to an optimal solution \cite{GleixnerSteffy2020}.

\paragraph{Failure due to stalling.} The second case occurs when the \ir method fails to significantly reduce violations for too many iterations in a row. By default, \ir aborts after two consecutive iterations where the maximum violation was decreased by a factor less than $2^4$.

\paragraph{Failure due to incorrect status.} The last two cases are when the \fp solver incorrectly detects an instance known to be feasible as infeasible, or an instance known to be bounded as unbounded.
Firstly, this situation can occur when the auxiliary LPs \eqref{eq:feasprob} for checking feasibility or \eqref{eq:unbdprob} for checking boundedness are solved exactly by \ir:
By construction, \eqref{eq:feasprob} is feasible and bounded, and \eqref{eq:unbdprob} is bounded.
Secondly, this can happen for the original LP after feasibility and boundedness has been established by exact solution of one or both of these auxiliary LPs.

\paragraph{}

In the original \ir method presented in \cite{GleixnerSteffyWolter2016}, running into numerical problems or incorrect status claims would enable a recovery mechanism that modifies certain settings such as presolving, scaling, tolerances, etc. Since this mechanism involves restarting the solving process from scratch, it can be very time-consuming and is disabled in our combined approach.
Instead, whenever one of these problems occurs, we boost the precision and restart the \ir procedure in higher precision.

\subsection{Boosting the precision}

\begin{figure}[t]
   \centering
\definecolor{DaLightBlue}{rgb}{0.63,0.79,0.95}
\definecolor{DaLightBrown}{rgb}{0.8,0.58,0.46}
\definecolor{DaLightGreen}{rgb}{0.56,0.93,0.56}
\definecolor{DaDarkBlue}{rgb}{0.0,0.0,0.46}
\begin{center}
   \begin{tikzpicture}[scale=.8]
      \tikzstyle{state} += [draw=black,rectangle, inner sep=1.5mm];
      \tikzstyle{link} += [->];

      \node (start) at (0.0,4.7) [state] {\scriptsize start};
      \node (fail) at (0.0,6.3) [state] {\scriptsize fail};
      \node (boost) at (5.5,6.3) [state,fill=DaLightBlue,rounded corners] {\scriptsize precision boosting step};
      \node (refine) at (5.5,4.2) [state, rounded corners] {\parbox{3.5em}{\scriptsize\centering iteratively\\refine LP}};
      \node (unb) at (3.1,-0.4) [state] {\scriptsize unbounded};
      \node (opt) at (5.5,-0.4) [state] {\scriptsize \smash[b]{optimal}};
      \node (inf) at (7.9,-0.4) [state, inner xsep=1.5ex] {\scriptsize infeasible};
      \node (stop) at (5.5,-2.3) [state, inner xsep=1ex, inner ysep=1.2ex] {\scriptsize \smash[b]{stop}};
      \node (testunbd) at (3.3,2.75) [state, rounded corners,fill=DaLightBrown] {\scriptsize test unboundedness};
      \node (testfeasunbd) at (3.3,1.15) [state, inner xsep=2.8ex,rounded corners,fill=DaLightGreen] {\scriptsize \smash[b]{test feasibility}};
      \node (testfeas) at (7.7,1.15) [fill=white, state, inner xsep=2.8ex, rounded corners,fill=DaLightGreen] {\scriptsize \smash[b]{test feasibility}};

      \draw [link] (start) -- ($(refine.west)+(0,0.5)$);
      \draw [link] (boost) -- node[anchor=south]{\tiny \hspace{3ex} limit reached} (fail);
      \draw [link] (boost) -- node[anchor=east]{\tiny boost successful} (refine);
      \draw [link] (refine) -- (opt); 
      \draw (unb.-63) edge [bend right=40] (stop.west);
      \draw (opt) -- (stop);
      \draw (inf.-117) edge [bend left=40] (stop.east);


      \draw [link] (refine.south west) -| node[anchor=south] {\tiny\hspace*{-3.5ex} fp unbounded} (testunbd.north);
      \draw [link] ($(refine.east)+(0,0.5)$)-| ++(1.6,0) -| node[above right=.15cm and .05cm, text width=1.5cm]{\tiny numerical troubles} ++(0,1) |- (boost.east);
      \draw [link] (testunbd.south) -- node[anchor=188] {\tiny unbounded} (testfeasunbd.north);
      \draw [link] (testunbd.232) -- node[anchor=east] {\tiny bounded} ++(0,-0.64) -| ++(-1.92,5.4)
      -| node[anchor=south, near start] {\tiny\hspace*{-2.9em} rejected unboundedness} (boost.125);
      \draw [link] (testfeasunbd.south) -- node[anchor=190.7] {\tiny feasible} (unb.63);
      \draw [link] (testfeasunbd.232) -- node[anchor=east] {\tiny infeasible\hspace*{0ex}} ++(0,-0.6) -| ++(-1.92,-3.35) -- ++(8.81,0) |- (inf.east);

      \draw [link] (refine.south east) -| node[anchor=south] {\tiny fp infeasible\hspace*{-1.5ex}} (testfeas.north);
      \draw [link] (testfeas.308) -- node[anchor=west] {\tiny feasible\hspace*{0ex}} ++(0,-0.6) -| ++(1.92,6.958)
      -| node[anchor=south, near start] {\tiny rejected infeasibility\hspace*{-0.95em}} (boost.55);
      \draw [link] (testfeas.south) -- node[anchor=-9] {\tiny infeasible} (inf.117);
   \end{tikzpicture}
\end{center}
   \caption{Flowchart illustrating the combined algorithm}
   \label{fig:generalflowshart}
\end{figure}
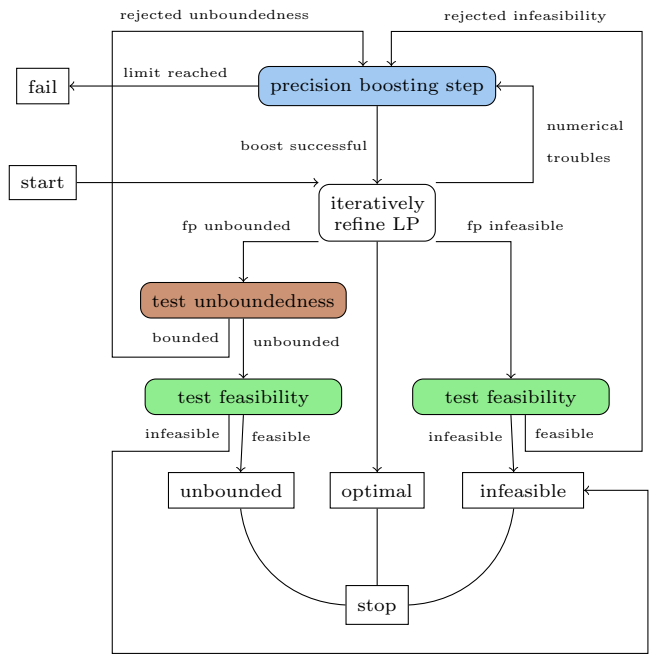

In the following, we explain important details of the combined algorithm. A flowchart illustrating the algorithm can be seen in Figure \ref{fig:generalflowshart}.

The precision boosting step can be split into three parts: First, the default arithmetic precision of all operations and data structures is increased to the new precision, then the LP is approximated in the increased precision from the rational LP, and finally the tolerances of the solver are decreased.

\paragraph{Increase precision.} We increase the precision very similarly to \qsoptex \cite{Espinoza2006}, which means first solving in double precision, then with $128,192,288, \ldots$~bits. Due to implementation details we currently do not support quad precision ($128$ bits), but instead directly increase to $192$~bits after the double-precision solve. A second difference is that we set the maximal precision limit at $1000$ bits, because above this precision some tolerances in \soplex{} (that are expressed in double precision) would automatically be rounded to zero. This is smaller than the $3164$~bit limit of \qsoptex, but note that the maximal precision is never reached in any of our experiments.

\paragraph{Load LP from rational LP.} Computing a more accurate approximation of the rational{} LP after each precision boost is necessary. In fact, in some instances, the roundoff-errors introduced when approximating the rational{} LP in \fp precision are the reason why the \fp solver returned a wrong basis. Consequently, after each \pb step, the rational coefficients $A, b, \ell, c$ of LP \eqref{eq:lpprimal} are approximated with the increasingly accurate \fp coefficients $\bar{A}, \bar{b}, \bar{\ell}, \bar{c}$.

\paragraph{Decrease tolerances.} The last part of the \pb step is the update of the tolerances in the \fp solver. Note that in contrast to a pure \pb approach, we are not required to change the feasibility and optimality tolerances of the \fp solver, since this is handled by the \ir procedure. However, the tolerances below which numbers are considered zero by the different components of the solver do need to be decreased.
Similarly to \qsoptex \cite{Espinoza2006}, we use constant factors to scale the tolerances in accordance with the current precision. Due to implementation details of the library we use to express higher-precision numbers, we express our tolerances in relation to the current precision, approximated in base $10$. For a precision of $p$ bits, let $p'$ be the such that $2^{-p} = 10^{-p'}$, \ie $p'=p \frac{log(2)}{log(10)}$. Then we set the tolerances as $10^{-\floor{p'*c}}$, where $c$ is a constant that can be different for each tolerance, \eg $c=1$ for the tolerance below which values are considered zero during the solve and $c=0.625$ for the tolerance below which values are considered as zero during pivot element selection.
The values of the constants are chosen to be consistent with the tolerance values in double precision.

After the precision has been increased, we restart the \ir procedure. It is clear that it is advantageous to not solve from scratch with higher-precision but rather to warm-start the solving process from an advanced basis.
If \pb was performed because of stalling, it is clear that we can restart from the last floating-point optimal basis. However, if \pb was performed because of numerical difficulties, it is desirable to restart with a more stable basis.
To achieve this, we store the basis at a geometrically increasing frequency: if the number of iterations is a power of two, and at least all \K{} iterations.

If unboundedness or infeasibility is detected by the \fp solver during the solution process, first one of the auxiliary LPs \eqref{eq:feasprob} or \eqref{eq:unbdprob} is created and solved exactly with the combined algorithm. If feasibility respectively boundedness has already been established by solving such an auxiliary problem, then any future claim of infeasibility respectively unboundedness immediately leads to a precision boosting step.
Furthermore, we store separate advanced bases for the original, the feasibility, and the unboundedness problem.

\subsection{Convergence guarantees of the combined algorithm}

The convergence proof of LP iterative refinement relies on the notion of a \emph{limited-precision LP oracle}~\cite{GleixnerSteffy2020}. An oracle is a limited-precision LP oracle, if there exist constants  $\sigma > 0$ and $0 < \eta < 1$, such that the oracle can, for any LP $\eqref{eq:lpprimal}$, produce an approximate primal-dual floating-point solution $(\bar x, \bar y)$ with
\begin{subequations}
   \label{eq:oracle-intext}
   \begin{align}
      \label{eq:oracle1-intext}
      \Norm{A\bar x - b}_{\infty} &\le \eta, \\
      \label{eq:oracle2-intext}
      \bar x &\ge \ell - \eta \onevec, \\
      \label{eq:oracle3-intext}
      c - A^T \bar y &\ge -\eta \onevec, \\
      \label{eq:oracle4-intext}
      \Abs{(\bar x - \ell)^T(c-A^T\bar y)} &\le \sigma.
   \end{align}
\end{subequations}
In other words, all primal and dual violations are bounded by $\eta$, and complementary slackness violations are bounded by $\sigma$.

Since our combined algorithm still uses \ir, it is clear that any convergence guarantee for the pure \ir approach immediately transfers to the new algorithm. If our base algorithm to solve the LP is restricted to double precision, then it is unreasonable to expect it to be able to act as a limited-precision LP oracle for all possible inputs.
Hence, \ir currently lacks a self-contained convergence guarantee.

This is where our combination with precision boosting does not only
provide a \emph{practical}, but also a \emph{theoretical} contribution.
If we allow to increase the precision to any arbitrary number, then it becomes possible to construct such an oracle and prove an unconditional convergence guarantee.

A detailed proof of convergence to a solution that satisfies $\eqref{eq:oracle1-intext}-\eqref{eq:oracle4-intext}$ is beyond the scope of this paper, but
%
intuitively it is clear that a solution with a primal and dual violation of at most $\eta < 1$ can be found, if the simplex can be executed in arbitrarily high precision and the largest numerical errors encountered during the entire course of the algorithm tend towards zero as the precision increases.
However, note that any such statement must depend on the implementation of the underlying numerical linear algebra routines, in particular the LU~update.
Many modern simplex implementations use Forrest-Tomlin updates~\cite{Forrest1972} in order to update the LU factorization of the basis matrix. While computationally very efficient due to its sparsity-preserving properties, this update is not proven to be backward stable. For a convergence proof, it would suffice to use the Bartels-Golub update~\cite{BartelsGolub1969} instead, which is known to be backward stable~\cite{Bartels1971}.

In this regard, note that to our knowledge the literature currently lacks a theoretical proof of convergence for pure precision boosting.
Such a proof would require special attention in decreasing the primal and dual feasibility tolerances. The current implementation of \qsoptex \cite{Espinoza2006}, \eg decreases the primal and dual feasibility tolerance by the same order of magnitude by which the precision is increased, \ie the relative difference between target tolerance and numerical accuracy remains unchanged.
This seems to work well in practice.
However, in order to prove convergence to an optimal basis theoretically, this relative difference would need to grow with every boosting step, \ie tolerances would need to decrease at a slower rate than the rate at which precision increases. The results in~\cite{Ogryczak1988} suggest that numerical stability is achieved as long as the tolerances exceed the forward error of the computed solution.

By contrast, our proof of convergence for \ir with precision boosting does not suffer from such interdependencies.
The reason is that precision boosting is only used to force the primal and dual violations below $\eta < 1$, and this $\eta$ can remain fixed.



\section{Computational Study}
\label{sec:comp}

We investigate the performance of our proposed combined algorithm, comparing it to both previously existing methods on a large set of LPs and MIPs. In particular, we are interested in the following questions. First, in the context of pure LP, \emph{how does the combined algorithm compare to pure \ir and pure \pb, and what are the strengths of the individual methods?} From the existing literature, we expect \ir to be faster on instances where it works, while \pb is supposed to handle numerically challenging problems with more consistency. As this is the first time both methods were implemented inside the same LP solver, we want to verify if this is still the case. Second, we want to determine on which sets of instances the respective algorithms perform better, and see if the combined algorithm can leverage the strengths of both parts. This is discussed in Section~\ref{sec:lpexp}.
 Third, in the context of MIP, we want to determine \emph{if the performance improvements from pure LP experiments translate to MIP solving, and if the combined algorithm can help to reduce the number of failed calls to the exact LP solver.} Especially the last part is of interest, since failing exact LPs was an issue when introducing numerical cutting planes to an exact MIP solver~\cite{EiflerGleixner2023}, where custom techniques were developed to make cuts numerically easier for the exact LP solver. This is discussed in Section~\ref{sec:mipexp}.

\subsection{Setup and test set}
The experiments were all performed on a cluster of Intel Xeon Gold 5122 CPUs with
3.6 GHz and 96 GB main memory. For all symbolic computations, we use
the GNU Multiple Precision Library (GMP) 6.1.2 \cite{GMP}.
All compared algorithms are implemented within \soplexv, and are freely available on GitHub\footnote{\url{https://github.com/scipopt/soplex}}. For exact MIP experiments we use a development version of \scip, which is also publically available\footnote{It can be obtained from \url{https://github.com/scipopt/scip/tree/exact-rational}.}~
 and which uses \papilov \cite{Papilo} for exact rational presolving.

For the pure LP tests, we use two different test sets. The first, which we call \irpaper, is a collection of instances from \cite{GleixnerSteffyWolter2016}, containing instances from the Netlib LP test set including the kennington folder, Hans Mittelmann's benchmark instances, Csaba M\'esz\'aros's LP collection, the LP relaxations of the COR@L
mixed-integer programming test set, and the LP relaxations of the MIPs from \miplib~ instances up to and including \miplib~2010.
The second test set, which we call \cutlib, is comprised of $100$ instances that all stem from subproblems encountered by exact \scip for instances of the \miplib~2017 benchmark set that proved difficult for the exact LP solver~\cite{EiflerGleixner2023}. 

For MIP experiments, we use the \miplib~2017 benchmark instances; in order to save computational effort, we exclude all those that could not be solved by the floating-point default version of \scipv within two hours. We use three random seeds for the remaining $132$ instances, making the size of our test set $396$. Note that for the pure LP experiments, we ran each instance only once since performance variability is not as pronounced in pure LP experiments. The time limit was set to $7200$ seconds for all experiments.
For the LP experiments we report aggregated times in shifted geometric mean with a shift of $0.1$ seconds, and LP iterations in shifted geometric mean with a shift of~$10$. For the MIP experiments, we use a shift of $1$ second for the time and $100$ nodes for the number of nodes.

\subsection{Pure LP experiments}
\label{sec:lpexp}

We compare the following three settings:
\begin{itemize}
   \item pure \ir with precision-boosting disabled, referred to as \pureir,
   \item pure precision boosting with \ir disabled, referred to as \purepb, and
   \item the proposed combination of both, referred to as \combi.
\end{itemize}
In all three settings, we use a rational factorization of the final basis matrix to test for exact optimality. Although there exist instances where the reconstruction approach presented in \cite{GleixnerSteffyWolter2016} performs better, in preliminary experiments the factorization approach clearly outperformed the reconstruction approach on average for all three settings.
Consequently, on instances where the final basis of the initial double-precision simplex solve is already optimal, there is no difference between the settings. We therefore only report results for instances where the final basis of the initial double-precision simplex solve is not optimal.  This leaves us with $84$ instances from \irpaper and $91$ instances from \cutlib. Furthermore, the time and iterations for the initial floating-point LP solve are the same for all three settings, since they all use the same floating-point simplex solver. We therefore only report the time and iterations \emph{after} the initial floating-point LP solve has finished.

As reported in Table~\ref{tbl:fpnotoptimal}, \combi as well as \purepb, were able to solve $79$~of the $84$ instances in \irpaper, with $5$ timeouts, while \pureir solved $44$ instances, and had to abort on the remaining $40$ instances. On \cutlib, \combi and \purepb were able to solve all $91$ instances, while \pureir solved $26$ instances. This demonstrates that some form of precision boosting is crucial for solving these numerically more difficult instances.

In terms of solving time, we compare all three settings on the subset of instances that were solved to optimality by all three settings, since the instances on which \pureir aborts unsuccessfully are not useful for comparison.

Table~\ref{tbl:alloptimal} reports aggregate results on both test sets. Comparing \combi and \pureir, we observe that \combi is slightly faster on \irpaper (by $\percentage{0.04}{1.15}$) and slower on \cutlib (by $\percentage{0.06}{0.37}$). On these instances, \purepb performs the worst, with a slowdown of $\increase{1.78}{1.15}$ on \irpaper, and $\increase{0.74}{0.37}$ on \cutlib, compared to \combi. Some of this can be explained by the selection of instances, which include the instances where \pureir successfully recovers from a faulty claim of infeasibility. More explicitly, the first floating-point solve claims infeasibility, but the feasibility test \eqref{eq:feasprob} shows that the instance is feasible. In that case, \pureir can sometimes recover by restarting, and those cases appear in the test set of all optimal instances, while the cases where the recovery is not successful are excluded. In that sense, looking at the subset of all optimal instances is biased in favor of \pureir. Nevertheless, it is clear that \pureir is faster than \purepb on instances that can be solved to optimality by both settings, while being comparable with \combi.

When the basis of the initial floating-point LP solve is not optimal, \purepb and \combi solve the same set of instances.
Hence, we can compare them fairly using the aggregate results from Table~\ref{tbl:fpnotoptimal}, which includes instances that fail with some setting.
On \irpaper, \combi is $\reduction{4.12}{4.98}$ faster than \pureir, while on \cutlib, \combi is $\increase{0.96}{0.91}$ slower than \pureir. This indicates that on sets of instances where the \ir algorithm already works well, it is beneficial to avoid precision boosting by running \ir first. On the other hand, in $71\%$ of the instances of \cutlib, \ir fails and the \combi algorithm needs to boost the precision. In these cases, simply performing precision boosting right from the start is often faster than trying \ir first.
However, if we split the union of both test sets into the subset of instances where \combi performed at least one precision boost, and the subset of instances where no precision boost was performed, we observe that \combi is $\reduction{3.03}{3.83}$ faster on the instances where no precision boost was performed, while only being $\reduction{6.52}{6.66}$ faster on the other subset.

\newcommand{\itin}{initial\xspace}
\newcommand{\itboost}{boosted\xspace}
\begin{table}
   \centering
   \caption{Comparison of the different variants on the subset of instances that where solved to optimality by all solvers, and where the final basis of the initial (double) precision LP solve is not optimal. Column ``\itin'' shows the number of simplex iterations in initial precision, column ``\itboost'' shows the number of iterations after the first precision boost.}
   \label{tbl:alloptimal}
   \begin{tabular}{llrrrrr}
        \toprule
        & & & & \multicolumn{2}{c}{iterations} \\
        \cmidrule{5-6}
        Test set  & setting & size & time & \itin & \itboost \\
        \midrule
        \irpaper & \combi  & 44   & 1.15 & 25.72 & -        \\
                 & \pureir & 44   & 1.19 & 25.72 & -        \\
                 & \purepb & 44   & 1.78 & 26.00 & 17.70    \\
        \midrule
        \cutlib  & \combi  & 26   & 0.37 & 26.91 & 0.71     \\
                 & \pureir & 26   & 0.31 & 26.57 & -        \\
                 & \purepb & 26   & 0.74 & 18.30 & 22.26    \\
        \bottomrule
\end{tabular}

\end{table}

\begin{table}
   \centering
   \caption{Comparison of the different variants on the subset of instances where the final basis of the initial (double) precision LP solve is not optimal. Column ``\itin'' shows the number of simplex iterations in initial precision, column ``\itboost'' shows the number of iterations after the first precision boost.}
   \label{tbl:fpnotoptimal}
   \begin{tabular}{llrrrrrr}
        \toprule
        & & & & & \multicolumn{2}{c}{iterations} \\
        \cmidrule{6-7}
        Test set  & setting & size & solved & time & \itin  & \itboost \\
        \midrule
        \irpaper & \combi  & 84   & 79     & 4.12 & 59.78  & 13.63    \\
                 & \pureir & 84   & 44     & 0.71 & 25.85  & -        \\
                 & \purepb & 84   & 79     & 5.98 & 52.19  & 44.66    \\
        \midrule
        \cutlib  & \combi  & 91   & 91     & 0.96 & 108.98 & 8.95     \\
                 & \pureir & 91   & 26     & 0.39 & 89.03  & -        \\
                 & \purepb & 91   & 91     & 0.91 & 37.75  & 36.79    \\
        \bottomrule
\end{tabular}


\end{table}

\begin{figure}
   \centering
   \caption{Comparison of solving times for \purepb vs. \combi. The left plot shows the results for \irpaper, the right plot shows the results for \cutlib. Instances where \combi performed at least one precision boost are marked in red.}
   \label{fig:speeds}
   \includegraphics[width=.49\textwidth]{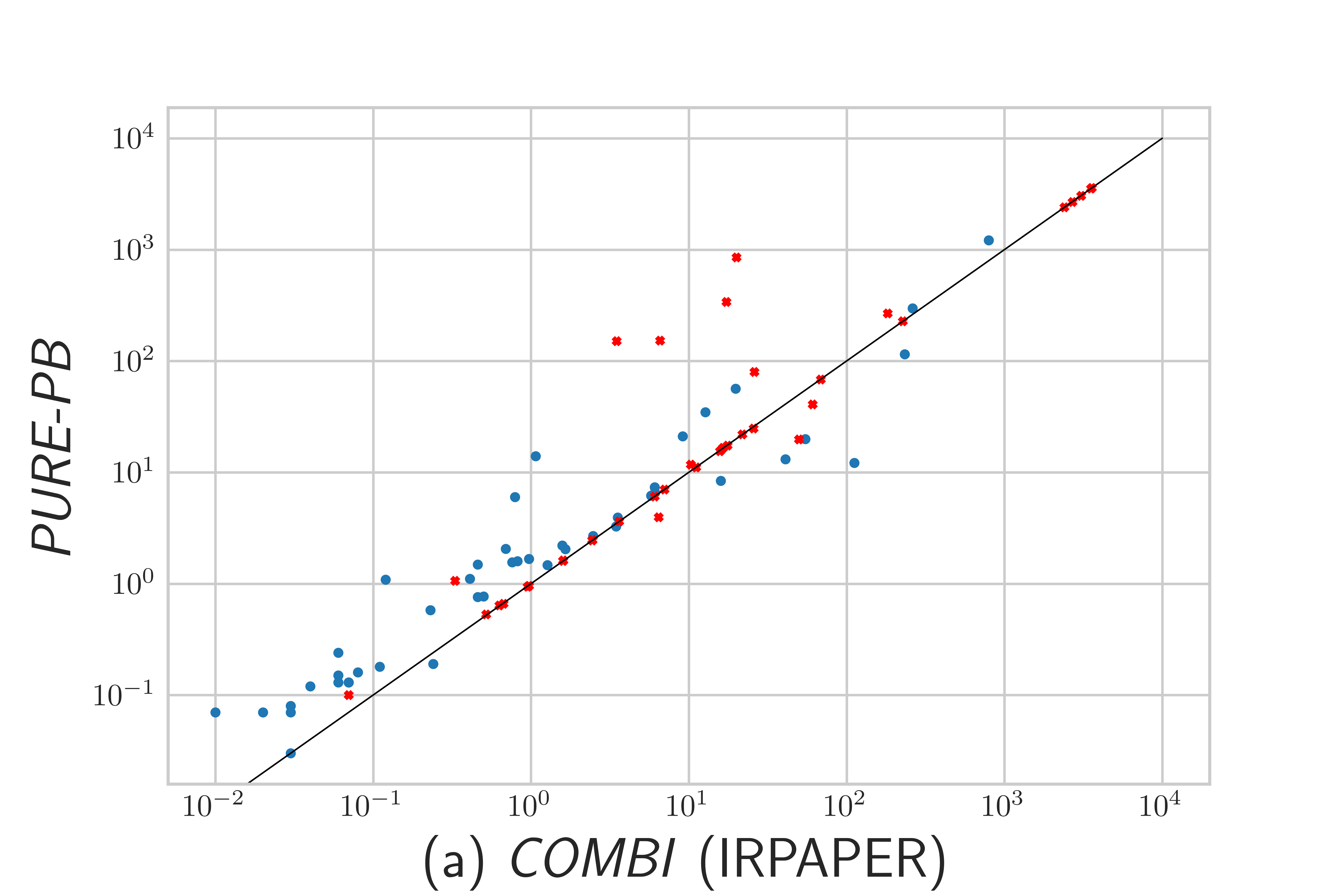}
   \includegraphics[width=.49\textwidth]{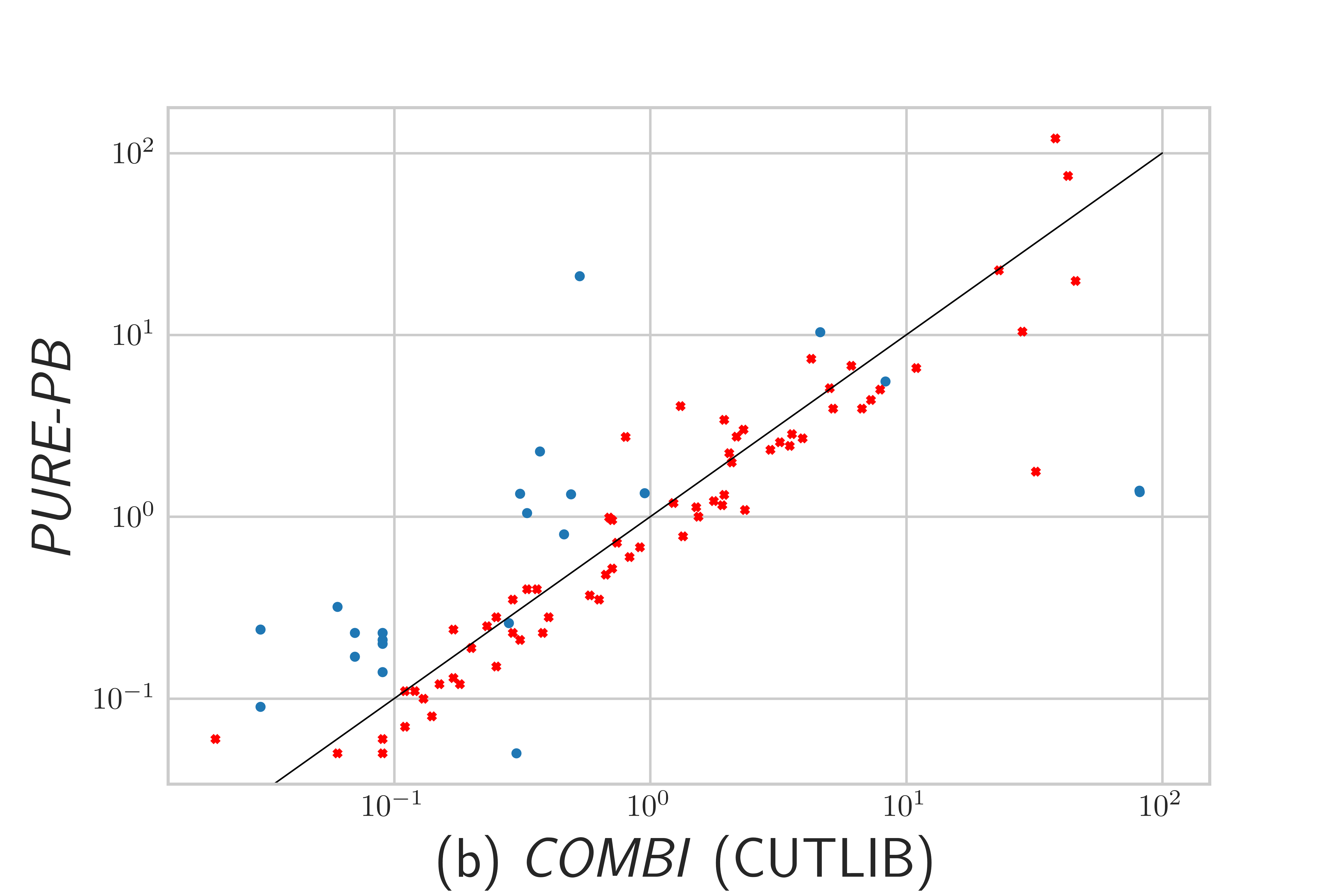}
\end{figure}

To summarize, overall \combi performs best among all settings, displaying a good tradeoff between speed on numerically easier instances and robustness on more difficult ones.
To further illustrate this point, Figure~\ref{fig:speeds} shows a scatter plot of solving times for \combi vs. \purepb. Points above the diagonal are instances where \combi is faster than \purepb, points below the diagonal are instances where \purepb is faster than \combi. Instances where \combi performed at least one precision boost are marked in red. Although there are outliers, it is clearly visible that \combi is faster than \purepb on the blue instances (where no precision boost was performed), while being slower on the red instances, especially on \cutlib (right plot).


\subsection{Experiments in an exact MIP solver}
\label{sec:mipexp}

Exact LP solving is an essential subroutine for solving MIPs exactly, certainly when following the bybrid-precision LP-based branch-and-bound approach proposed in~\cite{CookKochSteffyetal2013}, and whenever continuous variables are present.
Although our new proposed variant \combi performed well in pure LP experiments, it is not clear whether this will translate to exact MIP solving. On the one hand, the additional overhead of precision boosting might not be worthwhile in the context of MIP solving, where branching is an alternative to precision boosting. On the other hand, \combi may help to remove the need for custom-made cut-weakening techniques that were developed in \cite{EiflerGleixner2023} in order to reduce the encoding length of coefficents in the LP relaxation.

To investigate these question, we performed experiments using \scip on the \miplib~2017 benchmark instances with four different \scip settings:
\begin{itemize}
   \item the default best setting determined in \cite{EiflerGleixner2023}, precision boosting disabled, referred to as \enc
   \item the same setting, precision boosting enabled (in the \combi variant), referred to as \encboost
   \item cut-weakening disabled, precision boosting disabled, referred to as \safegmi
   \item cut-weakening disabled, precision boosting enabled, referred to as \safegmiboost
\end{itemize}

\begin{table}
   \centering
   \caption{Comparison of the four MIP settings on the \miplib~2017 benchmark instances. Only those instances are included for which at least one precision boost was performed for one of the settings. Column ``nboost'' shows the shifted geometric mean of the number of precision boosts performed by the respective setting, column ``exlpfails'' shows the shifted geometric mean of the number of times the exact LP solver failed to solve an exact LP.}
   \label{tbl:mip}

\begin{tabular}{lrrrrrrrrr}
   \toprule
   Setting       & size & solved & time    & exlptime & nboost & exlpfails \\
   \midrule
   \enc          & 146  & 44     & 4669.29 & 63.1     & 0.0    & 5.4       \\
   \encboost     & 146  & 45     & 4491.30 & 77.9     & 13.5   & 0.0       \\
   \safegmi      & 146  & 36     & 5153.33 & 208.9    & 0.0    & 50.8      \\
   \safegmiboost & 146  & 35     & 4924.84 & 204.7    & 73.0   & 0.0       \\
   \bottomrule
\end{tabular}

\end{table}

To reduce the amount of performance variability as much as possible, we restrict our analysis to the subset of $146$ instances where at least one precision-boosting step was performed by at least one of the settings. The results for all of these instances are reported in Table~\ref{tbl:mip}.

We observe that enabling precision boosting in the \enc setting leads to a speedup of $
\reduction{4491}{4669}$, with one more instance solved to optimality ($3$ gained, $2$ lost).
In the \safegmi setting it leads to a speedup of $\reduction{4924}{5153}$ with one less instance solved ($4$ lost, $3$ gained).
This very small difference in solving times is due to the large offset incurred from the large number of timeouts. For comparing solving times, Table~\ref{tbl:mip-oneopt} is more meaningful, since it only includes instances that were solved by at least one. Here, we observe a speedup of $\reduction{1757}{1973}$ for \encboost over \enc, and of $\reduction{2306}{2646}$ of \safegmiboost over \safegmi.

\begin{table}
   \centering
   \caption{Comparison of the four MIP settings on the \miplib~2017 benchmark instances. Only instances where at least one precision boost was performed for one of the settings, and where at least one settings managed to solve to optimality are included. }
   \label{tbl:mip-oneopt}
   \begin{tabular}{lrrrrr}
   \toprule
   Setting      & size & solved & time    & nodes   & exlptime \\
   \midrule
   \enc          & 49 & 44     & 1973.50 & 76975.1 & 26.9     \\
   \encboost     & 49 & 45     & 1757.81 & 72551.4 & 29.7     \\
   \safegmi      & 49 & 36     & 2646.87 & 82332.1 & 178.3    \\
   \safegmiboost & 49 & 35     & 2306.36 & 68470.1 & 153.8    \\
   \bottomrule
\end{tabular}

\end{table}

Finally, let us compare the performance of both settings that use precision boosting: \encboost and \safegmiboost. Although the number of failed exact LPs is zero in both versions, the number of necessary precision boosts is $\fraction{73}{13.5}$ times higher in \safegmiboost, leading to an increase of $\increase{204.7}{77.9}$ in exact LP time and of  $\increase{4924}{4491}$~in overall solving time, as reported in Table~\ref{tbl:mip}. We note that the large number of timeouts gives a strong offset to the average solving times given in this table. If we compare only on the subset of instances that could be solved to optimality by at least one setting, then the speedup of \encboost over \safegmiboost is $\increase{2306.36}{1757.81}$. This is comparable to the speedup of $\increase{2646.87}{1973.50}$ of \enc over \safegmi, which shows that precision boosting does not impact the importance of cut-weakening techniques, but is rather a complementary technique.


\section{Conclusion}
\label{sec:conc}

In this article, we presented an improvement to the \ir algorithm for solving LPs exactly over the rational numbers.
By integrating the precision boosting method inside an outer \ir loop, the combined algorithm is more robust on numerically challenging problems and does not suffer from any slowdown on numerically well-behaved instances. This addresses the major shortcoming of \ir for exact LP solving, namely the absence of a reliable recovery mechanism in case of numerical difficulties.

We analyze the performance of our new combined algorithm on two different LP testsets, comparing it with pure \ir, as well as pure precision boosting.
The results show that the combined algorithm outperforms the other methods for solving LPs exactly.  It solves more instances than a pure \ir approach, and is faster than pure precision boosting. Furthermore, we show that using the combined algorithm as a subroutine for solving exact MIPs, we are able to reduce the number of failed exact LP calls to zero, while remaining as fast as the pure \ir approach.

We see two directions for future improvement.
First, and most importantly, the rational LU factorization that is used to verify the exact optimality of a final basis can become a major bottleneck on some problem instances. A more involved roundoff error-free LU factorization algorithm, such as the one proposed by \cite{Lourenco2019}, could help to alleviate this problem.
Second, the missing intermediate step of increasing the precision to hardware-supported $128$ bit quad precision could speed up the combined algorithm, as it makes the first boosted iteration more efficient, especially since one precision boost is sufficient in most cases.

\section*{Acknowledgments}
We would like to thank Martin Sidaway and his supervisor Paolo Zuliani for their valuable insights into the theoretical convergence of pure precision boosting.



\bibliographystyle{plain}
\bibliography{status-report}

\end{document}